\documentclass{amsart}
\usepackage{amssymb}
\usepackage{graphicx}
\usepackage{enumerate}
\newtheorem{theorem}{Theorem}[section]
\newtheorem{lemma}[theorem]{Lemma}
\newtheorem{prop}[theorem]{Proposition}

\newtheorem{definition}[theorem]{Definition}
\newtheorem{proposition}[theorem]{Proposition}
\newtheorem{corollary}[theorem]{Corollary}  
\theoremstyle{definition}

\theoremstyle{remark}

\numberwithin{equation}{section}

\def\al{{\alpha}}
\def \th{{\theta}}

\def\P {{\bf P}}
\def\R{{\bf R}}
\def\C{{\bf C}}

\begin{document}

\title{On some analytical properties of stable densities}.

\author{Sonia Fourati}
\address{LMI-INSA Rouen + LPMA\\
76131 St Etienne du Rouvray\\ FRANCE
}
\email{soniafourati@insa-rouen.fr}

%    General info

\date{}

\begin{abstract}
L.Bondesson  \cite {bond} conjectured  that the  density of a
positive $\alpha$-stable  distribution is hyperbolically
completely monotone (HCM in short) if and only if
$\alpha\leq 1/2$. This   was
proved recently by  P. Bosch and Th. Simon, who also
conjectured a strengthened version of this result.
We disprove this conjecture as well as a
correlated conjecture of Bondesson, while giving a  
short  new proof of the initial  conjecture,   as a direct
consequence  of a new algebraic property of  HCM and
Generalized Gamma convolution densities (GGC in short)
which we establish.
\end{abstract}

\maketitle

{\bf R\'esum\'e}

L.Bondesson a conjectur\'e que la densit\'e d'une 
variable al\'eatoire $\alpha$-stable positive est hyperboliquement completement monotone (HCM)
si et seulement si $\alpha\leq 1/2$.  Ce r\'esultat  a \'et\'e r\'ecemment \'etabli par P.Bosh et Th.Simon qui ont aussi conjectur\'e une version plus forte de ce r\'esultat. 
Nous infirmons celle-ci ainsi qu'une autre  conjecture de L. Bondesson. 
Nous donnons aussi une courte et nouvelle preuve de la conjecture initiale, comme  cons\'equence directe d'une nouvelle propri\'et\'e alg\'ebrique  des 
fonctions HCM et des densit\'es gamma g\'en\'eralis\'ees (GGC) que nous \'etablissons.

 \section{Introduction}
This paper is concerned with the HCM property for stable distributions  and GGC random variables, whose definitions we recall below. Hyperbolically completely monotone functions (HCM in short)
 were  introduced by Lennart  Bondesson \cite{bond} in order to analyze infinitely divisible distributions. On the other hand, 
 the generalized gamma convolutions (GGC in short) introduced by O. Thorin \cite{thor}, are the  weak limits of finite convolutions of Gamma random variables. These notions are closely related, indeed the main example of  HCM functions are the Laplace transform  of  GGC  variables.  
L. Bondesson proved, in \cite{bond}, that the $\alpha$-stable positive random variables (denoted $S_{\alpha}$),  with density $g_{\alpha}$, are  GGC for all $\alpha\in ]0,1]$ and that
they have an HCM-density when $\alpha=n^{-1}$, for any integer $n\geq 2$. He also conjectured that this property holds for   all $\alpha\leq 1/2$. In a previous preprint \cite{Fou}  we proved   that the density (denoted $G_{\alpha}$)  of  $S_{\alpha}^{-\beta}$,  (with $\beta := {\alpha\over 1-\alpha}$) is  HCM if $\alpha\in [1/3, 1/2]$. This implies easily the HCM property of $g_{\alpha}$ for $\alpha $ in this range. Moreover,  it is easy to see that  $\beta$ is the largest real number for which this property  holds.
Recently,   Pierre Bosch and Thomas Simon \cite{Bosch}  proved the full  original Bondesson conjecture. Their proof makes use of the following result from  Bondesson \cite{bond2}  " The independent product or ratio of  two GGC random variables is again GGC" . Furthermore they  conjectured  that $G_{\alpha}$ is also an  HCM function   for   all $\alpha\leq 1/2$.
In the present paper, we prove that actually $G_{\alpha}$ is not HCM for $\alpha<1/3$.
Moreover, for  $\alpha\in ]1/2,1[$, using the fact that  ${e^{\delta x}\over G_{\alpha} (x)}$ is HCM (see also \cite{Fou}) ,  we obtain that $G_{\alpha}$ is not the density of a GGC random variable. Since   $g_{\alpha} $ is a GGC-density,
this  gives an example of  a GGC random variable $S_{\alpha}$ such that $S_{\alpha}^{\gamma}$ is  not GGC  and $|\gamma|>1$, thus  providing a negative answer to a question of   L. Bondesson \cite{bond2}. Finally, using Bondesson new remarkable property already mentioned, we prove  that the multiplicative convolution 
of an HCM function and a GGC density is again HCM.  As we show,   initial Bondesson's conjecture is an immediate consequence of this result.

The central result of this  paper  is a representation of $G_{\alpha}$ for all $\alpha\in ]0,1]$. One consequence of this representation is the estimate of this density by a convex combination of   two gamma densities, namely $\Gamma (1/2, \delta)$ and $\Gamma (\alpha,\delta)$, with $\delta = (1-\alpha)\alpha^{\alpha\over 1-\alpha}$. Hopefully, this might be useful   for the  numerical investigation of these functions.

This paper is organized as follow :  In the first  part, we recall the facts on HCM functions and GGC random variables which are used in the sequel, we refer to Bondesson   \cite{bond} and James-Roynette-Yor \cite {Yor-Roy-James} for more general informations on the subject, other examples of GGC or HCM functions can be found in the recent work of W. Jedidi and Th. Simon \cite{WissSim} . In the second part, we do the same for $\alpha$-stable densities, for which we refer to Zolotarev \cite{zol}.

In part 3, we give a first rough  estimate of $G_{\alpha}$, as a function of a complex variable, obtained using the saddle point method.
In part 4, the main result is given : we give  a representation of the density $G_{\alpha}$ of $S_{\alpha}^{-\beta}$  for every $\alpha\in ]0,1[$ and,
as  a first corollary,  an estimate of $S_{\alpha}^{-\beta}$   by a convex combination of two Gamma distributions is given. 
In part 5, we prove that  $G_{\alpha}$ is not HCM for  $\alpha<1/3$ (which disproves the Simon-Bosch conjecture), is  HCM  for $\alpha\in [1/3, 1/2]$ and anti-HCM (see definition below) for $\alpha>1/2$.
As  consequences, more  precise estimates are given  when  $\alpha\in [1/3,1/2]$ and a corollary of the anti-HCM property for $\alpha>1/2$ is that  $S_{\alpha}$ is GGC although $S_{\alpha}^{-\beta}$ is not. 
In part 6, we prove that the convolution product (sometimes called mixing) of an HCM function and a GGC density is again HCM.  As a consequence, we obtain a short proof of   the HCM property of $S_{\alpha}$.

\section{Preliminaries}
\subsection{Hyperbolically completely monotone functions}\label{HCM}
We recall here the basic definition and properties of the class of  hyperbolically completely monotone functions, and refer to  \cite{bond} for more details.

\begin{definition}
A real positive valued function $H$ defined on $]0,+\infty[$ is called hyperbolically completely monotone (HCM) if, for every $u>0$ the function $H(uv)H(uv^{-1})$ is a completely monotone function of the variable $v+v^{-1}$.
\end{definition}
 
 Bondesson \cite{bond}  has obtained the following characterization  of HCM functions.

\begin{prop}\label{PropositionHCM}\ \par
$H$ is HCM if and only if it admits the following representation
\begin{equation}\label{HCM}
H(x)=cx^{\beta-1}\exp\left(-ax-\int_1^\infty\log\frac{x+t}{1+t}\mu_1(dt)-bx^{-1}-
\int_1^\infty\log\frac{x^{-1}+t}{1+t}\mu_2(dt)\right)
\end{equation}
where $ a, b,c$ are non negative constants and $\mu_1$ and $\mu_2$ are positive Radon measures on $[1,+\infty[$ that integrate $1/t$ at infinity.
\end{prop}
We shall use a slightly different but  equivalent representation of HCM functions,  obtained by an integration by part from (\ref{PropositionHCM}).  Denote $$\th(t)=\mu_1(]1,t]) 1_{t\geq 1}- \mu_2(]1, 1/t[)1_{t<1} +(\beta-1) $$
then
 $\th$ is a
(signed)  non decreasing c\`adl\`ag function. The following is an immediate consequence of Proposition \ref{PropositionHCM}.
\begin{corollary} $H$ is HCM if and only if it admits the following representation,
\begin{equation} \label{HCMb}
H(x) = c\exp ( - a x -b x^{-1})\exp\left( \int_0^{\infty} ({1\over x+t}-{1\over t+1}) \theta(t)dt \right)
\end{equation}

where $ a, b,c$ are non negative constants and $\th$ is a signed non decreasing function such that  $\int_0^{\infty} (1\wedge{1\over t^2})|\th(t)| dt <+\infty$ 

\end{corollary}

One has
$H(1) =ce^{-a-b}$. Moreover,  if $\th= \th_0$ is  a constant function and if  $a=b=0$,  then 
 $H(x)=cx^{-\th_0}$.   The integral condition $\int_0^{\infty} (1\wedge{1\over t^2})|\th(t)| dt <+\infty$ is the minimal condition  to ensure   finite values for $H(x)$ for every $x>0$. Note also that $H(x)$ may  be infinite   at $x=0$ and $x=+\infty$.
In the sequel,  the functions admitting this representation with $\th$ a non increasing function instead of a non decreasing function will be  called {\bf anti-HCM functions}.

 The representation (\ref{HCMb}) implies that $H$ has an analytic continuation on $\C\backslash ]-\infty,0]$. If we denote this continuation by $H$ again one has, using well known properties of the Stieltjes-Cauchy tranform:  
  $$H(-r^+):= \lim_{z\to -r, \Im(z)>0} H(z)= R(r) e^{-i\pi\th(r)}$$ 
 
 where $R(r) e^{-i\pi\th(r)}$ is the polar decomposition of the complex number $H(-r^+)$.  This  property will play a crucial role in the sequel.

 \begin{definition} 
 The {\sl generalized Gamma convolutions } (GGC) are the  random variables which belong to the smallest class containing  Gamma distributions and closed under taking sums of independent variables and weak limits.
 \end{definition}
 The following results can be found in Bondesson \cite{bond}

\begin{proposition}\label{HCMGGC}\

\begin{enumerate}\item
 A random variable  is GGC if and only if its Laplace transform is  an HCM function. 
 \item
  An  HCM function $H$  is the Laplace transform of a random variable  if and only $b=0,H(0)=1$,  and  the function $\th$ in the representation (\ref{HCMb})  is non negative. 
 Moreover, when these properties  are satisfied, $H$ is the Laplace transform of a GGC random variable.
 \end{enumerate}
 \end{proposition}
 
The class of HCM functions  and GGC random variables,  have been  much studied. We refer mainly to   Bondesson monography \cite{bond} and to Yor-Roynette-James \cite{Yor-Roy-James}, for GGC- random variables. 
Note that $e^{-x^{\alpha}}$, $\alpha<1$, is  an HCM function with $a=b=0$ and $\th(t)= \sin\pi\alpha t^{\alpha}$, while $e^{-x^{\alpha}}$  is the Laplace transform of the positive $\alpha$-stable distribution $S_{\alpha}$, thus $S_{\alpha}$ is GGC.

\subsection{Stable random variables}

Let $\alpha\in]0,2[,\rho\in]0,1]$, and suppose $\gamma=\alpha\rho\in [0,1]$ and let $g_{\alpha,\gamma}(x)$  denote  the density of the normalized   $\alpha$-stable random variable $S_{\alpha,\rho}$
 with  asymmetry parameter $\rho$ ($\rho= \P(S_{(\alpha,\rho)}>0)$ ) cf \cite{zol}). For $\rho=1$ and $\gamma=\alpha\rho=\alpha\in ]0,1[$ (and only for these values) this distribution is supported on the half axis $]0,+\infty[$ and we 
simply put $g_\alpha=g_{\alpha,\alpha}$. 

The function $g_{\alpha,\gamma}$ has Fourier transform 

$$e^{-(iu)^{\gamma} (-iu)^{\alpha-\gamma}}= \int_{-\infty}^{+\infty} e^{-iut}g_{\alpha,\gamma}(t) dt$$

where  $t^{\alpha}=\exp(\alpha\log(t))$ with $\log$ the principal determination of the logarithm.

The following  integral representation  (cf Zolotarev \cite{zol})  can be easily obtained  by Fourier inversion. 

\begin{lemma}{\bf Zolotarev}  For $r>0, 0\leq \alpha<1$ and $0\leq\gamma\leq 1$

\begin{equation}\label{Zolotarev}
g_{\alpha , \gamma}(r)=(2 i \pi )^{-1} \int_0^{\infty} (e^{-rt- t^{\alpha} e^{i\pi\gamma}}- e^{-rt- t^{\alpha}e^{- i\pi\gamma} })dt \end{equation}

\end{lemma} 

The above integral is  well defined for all $\alpha \in ]0, 1[$ and $|\gamma|\leq 1$. We will use it as   a definition in these cases.

  \begin{lemma}The function
$\tilde g_{\alpha,\gamma}(x)=x^{-1-\alpha}g_{\alpha,\gamma}(x^{-1})$ is 
\begin{enumerate}\item
decreasing on $]0,+\infty[$ if $0\leq \gamma \leq \alpha\leq 1 $.
\item
   completely monotone if $0\leq \gamma \leq \alpha\leq 1/2$
 \end{enumerate}
\end{lemma}

\begin{proof} Recall that, if $X$ is a   stable variable  with parameters $(2\alpha,\rho)$ and $Y$  is  an independent stable variable with parameters $(1/2,1)$, then $Z=XY^{\frac{1}{2\alpha}}$ is a stable variable with parameters $(\alpha,\rho)$ . Since 
$g_{1/2}(t)=\frac{e^{-\frac{1}{2t}}}{\sqrt{2\pi t^3}}$ one has
$$g_{\alpha,\gamma}(x)=2\alpha\int_0^\infty 
g_{2\alpha,\gamma}(y)
\frac
{e^{-\frac{1}{2}(y/x)^{2\alpha}}y^{\alpha}}
{\sqrt{2\pi}x^{\alpha+1}}
dy$$ Therefore
$$x^{-1-\alpha}g_{\alpha,\rho}(x^{-1})=2\alpha\int_0^\infty g_{2\alpha,\rho}(y)
\frac{e^{-\frac{1}{2}(yx)^{2\alpha}}y^{\alpha}}{\sqrt{2\pi}}dy$$
which is  decreasing in $x$ and completely monote if $2\alpha\leq 1$. 
\end{proof}

\begin{lemma} For $\alpha\leq \delta<1$ 

$$\int_0^{\infty}  g_{\alpha,\gamma}(xy)g_{\delta }(y)ydy =x^{\delta-1} g_{{\alpha\over\delta} , \gamma} (x^{\delta})  $$

\end{lemma}
\begin{proof} Let $\overline g$ denote the tail function of $g$,
$$\overline g(x)=\int_x^\infty g(y)dy.$$

Instead of the identity of the lemma, we rather prove the equivalent identity on the associated tail functions, 

$$\int_0^{\infty} \overline g_{\alpha,\gamma}(xy)g_{\delta }(y)dy =\frac{1}{\delta}\overline g_{{\alpha\over\delta} , \gamma} (x^{\delta})  $$

By (\ref{Zolotarev}) 

$$\overline g_{\alpha , \gamma}(r)=(2 i \pi )^{-1} \int_0^{\infty} (e^{-rt- t^{\alpha} e^{i\pi\gamma}}- e^{-rt- t^{\alpha}e^{- i\pi\gamma} })\frac{dt}{t}$$
therefore
$$\int_0^{+\infty} \overline g_{\alpha,\gamma}(xy)g_{\delta }(y)dy=(2 i \pi )^{-1} \int_0^{\infty} (e^{-x^\delta t^\delta- t^{\alpha} e^{i\pi\gamma}}- e^{-x^\delta t^\delta- t^{\alpha}e^{- i\pi\gamma} })\frac{dt}{t}$$
The proof follows by a simple change of variable ($t\to t^{1/\delta}$) in the integral and (\ref{Zolotarev}) again.
\end{proof}

\begin{proposition}\label {Monot} 
If $\alpha<\gamma\wedge 1/2< 1$ then $g_{\alpha,\gamma} (x)$ is not of constant sign. 
\end{proposition} 

\begin{proof}

If $\gamma\leq 1/2$, let $\delta={\alpha\over \gamma}< 1$ then  the complete monotonicity of $y^{-1-\delta}g_{\delta }(y^{-1})$   and the positivity of $g_{\alpha,\gamma} (x)$
 would imply the complete monotonicity of $g_{\gamma, \gamma}(x^{\gamma})$, but  this cannot be  true since $g_{\gamma, \gamma}(x^{\gamma})$ is not monotonous.
  
If $\gamma >1/2$ and $\alpha\leq  1/2$  then take $\delta=\alpha$ and obtain that $g_{1, \gamma} (x^{\alpha})$ would be completely monotonous and this is not true.  Actually this function is not  monotonous.

 \end{proof}
{\bf Remark}  Iterating the convolution with $\delta=1/2$ sufficiently, we could obtain that $g_{\alpha,\gamma} (x)$ is not of constant sign, for all $\alpha<\gamma\leq 1$.

Let $\alpha\in ]0,1[$,
$\beta:= {\alpha\over 1-\alpha}$
  and,  for all $x\in ]0,+\infty[$, 
  $$G_{\alpha}(x) :=  \beta^{-1}x^{-{1\over \alpha} } g_{\alpha} (x^{-\beta^{-1}})$$

  The  function $G_{\alpha}$ is the density of the distribution of $S_{\alpha}^{-\beta}$. It will play an important role in this paper.  The following integral representation
 \begin{equation}  \label {Galt}
G_{\alpha}(x)=(2i\pi \beta x)^{-1} \int_0^{\infty} (e^ {- t- e^{i\pi\alpha} t^{\alpha}x^{1-\alpha}} -  e^{- t- e^{-i\pi\alpha} t^{\alpha}x^{1-\alpha} } )dt\qquad x>0
   \end{equation}
  shows that $G_{\alpha}$
     has an analytic continuation to $C\backslash ]-\infty, 0]$, still denoted  $G_{\alpha}(z)$.

\section{A rough estimate of $G_{\alpha}$}  

Let $t_0$  be the the minimum of the function $ f(t)= t-t^{\alpha}$ for $t\in ]0,+\infty[$, and $\delta=-f(t_0)$, i.e. $t_0= \alpha^{1\over 1-\alpha}$ and $\delta=(1-\alpha)\alpha^{\alpha\over 1-\alpha}$.
Define $f_0(t)=f(t)-f(t_0)$.
The next  lemma gives a technical intermediate result that  will be improved in the next section.

\begin{lemma} \label{RghE} There exist constants $A>0$ and $B>0$ such that for all $z\in \C\backslash]-\infty,0]$

$$|G_{\alpha} (z) e^{\delta z}|\leq A+B |z|^{-1} $$

\end{lemma} 
Before proving this Lemma we need a new representation of $G_\alpha$. 
Consider the analytic function on ${\bf C}\setminus[0,+\infty[$ which coincides with the principal determination $t^{\al}$ on the upper half plane.
Let $f_+(t)$ be  the function obtained from $f_0$ by replacing  $t^{\al}$ by this function. In other words, $$f_+(t)= f_0(t)= t-t^{\alpha} +\delta \quad \hbox{  if }\quad \Im(t)>0$$

$$f_+(t)=t-e^{2i\pi \alpha} t^{\alpha} +\delta\quad \hbox{  if }\quad \Im(t)<0.$$

Similarly, let

$$f_-(t)=t-e^{-2i\pi \alpha} t^{\alpha} +\delta\quad \hbox{  if }\quad \Im(t)>0$$ 

$$f_-(t)= f_0(t)= t-t^{\alpha} +\delta \quad \hbox{  if }\quad \Im(t)<0.$$
One has $$f_+(\overline z)=\overline{f_-(z)}$$

\begin{lemma}
Let $\th  \in [-1, -(1/2-\alpha)^+[\cup  ](1/2-\alpha)^+,1]$, there exist two continuous  complex valued functions, $v^+_{\th}(r)$, $v^-_{\th}(r)$,  defined for $r\in [0,+\infty[$,  
such that $v^+_{\th}(0)=v^-_{\th}(0)=t_0$ and for all $r>0$, $v^+_{\th}(r)\not=v^-_{\th}(r)$ and 
$$f_+(v^+_{\th} (r)) = f_-(v^-_{\th} (r)) = re^{i\pi\th}$$
\end{lemma}

\begin{proof} 
Fix $\th\in ] 0,1]$, we shall build $v^+_{\th}(r)$ and $v^-_{\th}(r)$.
The point $t_0$ is a non degenerate saddle point for $f_0$ and $f_0(t_0)=0$, and $f_0(t) \sim {\alpha (1-\alpha) \over 2}  (t-t_0)^2$ in a neighborhood of $t_0$. By the  implicit function theorem there exists two distinct solutions $z$ satisfying the equation 
$$f_0(z)= re^{i\pi\th}$$ for $r$ small enough, moreover,  one of the two solutions, $v^+_{\th} (r)$ is in the upper half plane  and the other one, $v^-_{\th} (r) $, is in the lower half  plane.
Since $f_+=f_0$ on the upper half plane and $f_-=f_0$ on the lower half  plane, the complex numbers $v^+_{\th}(r)$ and $v^-_{\th} (r)$ satisfy the equations $$f_0(v^+_{\th} (r))= f_+(v^+_{\th} (r))= re^{i\pi \th}\qquad f_0(v^-_{\th} (r))=f_-(v^-_{\th} (r)))= re^{i\pi \th}$$
 Let $\th\in ](1/2-\alpha)^+,1]$. We  prove  that the two functions $v^+_{\th}(r)$ and $v^-_{\th} (r)$  can be continued for all $r>0$ and  $\th\in ](1/2-\alpha)^+,1]$. The case 
$-\th\in ](1/2-\alpha)^+,1]$ is obtained by conjugation. 
Since $f_+$  is analytic on the upper half plane,  by the open mapping  theorem, the function $v^+_{\th}$ can be continued as long as $v^+_{\th}(r) $ does not meet the real line.
The boundary values of the function $f_0$ on the upper half plane have imaginary part $0$ on the positive real line and negative imaginary part  on the negative real line, therefore  $v^+_{\th} (r)$ cannot converge to a real point since $\Im(f_+(v^+_{\th} (r)))= r\sin\pi\th>0$.

Similarly  $v^-_{\th}$ can be continued as long as $v^-_{\th}(r)$  does not reach the cut $[0,+\infty[$.
Let $t\in \R^+$, the boundary values of $f_-$ at $t$ are

\begin{eqnarray*}f_-(t^-)&=&f_0(t)= t-t^{\alpha}+\delta \in \R^+\\ 
f_-(t^+)&=& t- t^{\alpha} e^{-2i\pi \alpha} +\delta =( t-t^{\alpha} +\delta) + (1-e^{-2i\pi\alpha}) t^{\alpha}\\
&= &(t- t^{\alpha} +\delta) + 2 t^{\alpha}\sin\pi\alpha e^{i\pi (1/2-\alpha) } \end{eqnarray*}

Since $t-t^{\alpha} +\delta \in {\bf R}^+$, we see that $f_-(t^+)$ and $f_-(t^-)$  are in the cone $\{z; |{\arg} z| \leq |1/2-\alpha|\}$. On the other hand $f_-(v^-_{\th} (r))=re^{i\pi\th}$ always remains outside this cone, so that $v^-_{\th} (r)$ is defined for all $r>0$.
 \end{proof}

\begin{lemma} \label{Estv}\

\begin{enumerate}\item 
 There exists a  positive constant $A$  such that for $r\in \R_+$ and $\th \not\in  [ -(1/2-\alpha)^+ , (1/2-\alpha)^+]$,
  $$|v^+_{\th}(r)|\leq A+ 2 r\qquad  |v^-_{\th}(r)|\leq A+2r$$  
\item  
   $v^+_{\th}(r)\sim re^{i\pi\th}$ and $v^+_{\th}(r)\sim re^{i\pi\th}$  for $r\to +\infty$ and  $|\th|>(1/2-\alpha)^+$.
 \end{enumerate}
 \end{lemma} 
 \begin{proof} 
 This follows easily from the fact that there exists a constant $C>0$ such that
  $|{f_\pm(z)\over z}-1|\leq C|z|^{\alpha-1}$.
\end{proof}

  \begin{lemma}  For  $\th\in ^c [-({1/2-\alpha})^+,(1/2-\alpha)^-]$ and $\Re(ze^{i\pi\th})<0$,
    $$G_{\alpha} (z) =(2i\pi\beta )^{-1}e^{-\delta z}ze^{i\pi \th}.\int_0^{\infty} e^{z t e^{i\pi \th}} (v^+_{\th}(t)-v^-_{\th}(t)) dt$$
 \end{lemma} 

 \begin{proof} 
For all $u<0$ and $z\not\in ]-\infty,0]$ one has  
 \begin{eqnarray*}u - u^{\alpha} e^{i\pi\alpha}z^{1-\alpha}+\delta z&=& zf_+({u\over z})\\
u - u^{\alpha} e^{-i\pi\alpha}z^{1-\alpha} +\delta&=& zf_-({u\over z})\end{eqnarray*}
 Using this, we obtain from (\ref{Galt}) 
 $$G_{\alpha}(z)=(2i\pi\beta)^{-1}(\int_{1/2D_{\hat h }} e^{zf_+(t)} - \int_{1/2D_{\hat h } } e^{zf_-(t)} dt) $$
 where $1/2D_{\hat  h}$ is the half line $ \{-te^{i\pi h}; t\in[0,\infty[ \}$, $h$ is the argument of $z$. 
 We change again the contour and replace  the half line $1/2 D_{\hat h}$ by  the curve $[0,t_0]^+\cup \{v^+_{h} (s), s=0\to +\infty[\}$ for the first integral and $t\to [0, t_0]^{-}\cup \{v_{h} (s); ]-\infty, s]\}$ in the second one. Notice  also that  $f_+(t^+)=f_-(t^-)=f_0(t)$ for all $t\in [0,t_0]$, consequently, the  contribution of the two integrals over $[0,t_0]$ compensate each other and the end point of the half line and the two curves coincide at infinity.
 Finally  we obtain by the use of Cauchy theory that 
$$G_{\alpha}(z)=(2i\pi\beta)^{-1}(\int_0^{+\infty} e^{zf_+(v^+_{\th}(s))}dv^+_{\th}(s)- \int_0^{+\infty}  e^{zf_-(v^-_{\th}(s))} dv^+_{\th}(s)) $$
$$=(2i\pi\beta)^{-1}(\int_0^{+\infty}  e^{zse^{i\pi\th}}(dv_{\th}(s)- dv_{\th}(s))$$
 The integral representation of $G_{\alpha}$ follows after an integration by part.

\end{proof}

 {\sl  Proof of Lemma \ref{RghE}}.

   If $z\not\in [0,+\infty[$ let  $h$ be such that $ze^{i\pi h}= -|z| e^{i\pi \varepsilon}$  with 
   $|\varepsilon|\leq |1/2-\alpha| $.  One has 
   
     $$e^{\delta z}G_{\alpha} (z) =(2i\pi \beta)^{-1}z. \int_0^{+\infty} e^{-|z| t.e^{i\pi h} } (v^+_{\th}(t)-v^-_{\th} (t)) dt$$
  Using the estimate of $v^+_{\th}(t)$ and $v^-_{\th} (t)$ given in  lemma \ref{Estv}  we obtain, 

 $$|e^{\delta z}G_{\alpha} (z) |\leq (2\pi\beta )^{-1} (A (\sin\pi\alpha)^{-1} + 2(\sin\pi \alpha)^{-2} |z|^{-1} )$$
\qed

\section{The main result}

\begin{theorem} \label{main}
 There exists a continuous   function $\th$, taking values in $]0,1[$, such that,  
 for all $z\in{\bf C}\setminus ]-\infty,0]$, $\alpha\in ]0,1[$, 
 
\begin{equation}
G_\alpha(z)=G_{\alpha}(1)  e^{-\delta (z-1)} \exp \int_0^{\infty} ({1\over z+t}- {1\over t+1})\th(t) dt 
\end{equation}
Moreover, 
\begin{eqnarray*}
G_{\alpha}(z)\sim  c_0 z^{-\alpha} (1+O(z^{1-\alpha})\qquad z\to 0\\
G_{\alpha}(z)\sim c_{\infty} z^{-1/2} e^{-\delta z}  (1 +O(z^{-1})) \qquad z\to\infty 
\end{eqnarray*}
$$\text{with}\quad c_{\infty}= (2\pi \beta) ^{-1/2} \alpha^{\beta/2}. \qquad  c_0=(2\pi\beta)^{-1}\Gamma (\alpha+1)\sin\pi\alpha $$
\end{theorem} 

 The following estimate of $G_{\alpha}(x)$ on the real line is an immediate consequence of  this representation. 

\begin{corollary}\label{sim} Let \begin{eqnarray*}A_{\pm}& =&\sup_{x\in [0,1]} [x^{\alpha} e^{\delta x }G_{\alpha}(x)]^{\pm 1}\\ 
B_{\pm}& =&\sup_{x\in ]1,+\infty[} [x^{1/2} e^{\delta x }G_{\alpha}(x)]^{\pm 1}\end{eqnarray*}
then   $A_+,B_+,A_-, B_-$ are finite and non zero. 
Moreover, let $$f_1(x)=x^{-\alpha}e^{-\delta x}1_{[0,1]}(x), \quad f_2(x)=x^{-1/2}e^{-\delta x} 1_{]1,+\infty[} $$
then 
$$A_- f_1(x) + B_-f_2(x) \leq G_{\alpha}(x)\leq A_+ f_1(x) + B_+f_2(x) $$
\end{corollary}

For the  proof of Theorem \ref{main} we need first to study the behavior of $G_{\alpha}$ near the boundary $]-\infty,0[$.  Using 
(\ref {Galt}) one gets

\begin{equation}\label{Galt+}G_{\alpha}(-r^+) =(2i\pi\beta)^{-1}\int_0^{\infty} e^{-rt} (e^{rt^{\alpha}} - e^{rt^{\alpha} e^{-2i\pi\alpha}})dt \end{equation}

\begin{prop}\label{Galpha} For $r>0$ one has
\begin{enumerate}
\item  $\Im(G_\alpha(-r^+))<0$.

\item $G_{\alpha} (-r^+)= c_0 r^{-\alpha} e^{-i\pi \alpha} ( 1+ O(r^{1-\alpha}) )$ for $r\to 0$

\item $G_{\alpha} (-r^+) =-ic_{\infty}r^{-1/2}e^{\delta r} (1+ O(r^{-1}) )$ for $r\to \infty$.
\end{enumerate}
\end{prop}

\begin{proof} 

 -{\sl (1)}   follows from  $$\Re (e^{rt^{\alpha}} - e^{rt^{\alpha} e^{-2i\pi\alpha}})=e^{rt^{\alpha}} - e^{rt^{\alpha} \cos 2\piÊ\alpha}\cos [r t^{\alpha} \sin(2\pi\alpha)]\geq 0$$

 -{\sl (2)} The change of variables $t\to {t\over r}$ in (\ref{Galt+}) gives
 
 $$G_{\alpha} ( -r^+)= (2i\pi \beta r)^{-1}  \int_0^{+\infty}  e^{r^{1-\alpha} t^{\alpha} } - e^{r^{1-\alpha} t^{\alpha}e^{-2i\pi\alpha}}e^{-t}dt$$
The function
$$E(z)=(2i\pi \beta z)^{-1}  \int_0^{+\infty}  e^{z t^{\alpha} }- e^{z t^{\alpha}e^{-2i\pi\alpha}}e^{-t}dt$$

is entire and $E(0)=(\pi\beta)^{-1}e^{-i\pi\alpha}\Gamma(\alpha+1)\sin\pi\alpha$, moreover one has
 
 $G_{\alpha} (-r^+)= r^{-\alpha} E(r^{1-\alpha})$ from which  {\sl (2)}  follows.

-{\sl (3)} Using  Laplace method  we  obtain the following estimate   

$$ (2i\pi \beta)^{-1}\int_0^{\infty} e^{-r(t-t^{\alpha})}dt =-ic_{\infty}r^{-1/2} e^{\delta r} (1+O(r^{-1}))\qquad r\to\infty$$ 

Moreover,  
$$|e^{-\delta r}G_{\alpha}(-r^+)-(2i\pi\beta)^{-1}\int_0^{\infty} e^{-r(t-t^{\alpha}+\delta)}dt|\leq  (2\pi\beta)^{-1}\int_0^{\infty}e^{-r(t-t^{\alpha} +\delta)} e^{ -(1-\cos 2\pi\alpha) t^{\alpha}r }dt. $$ Since 
$e^{-r(t-t^{\alpha} +\delta)}\leq 1$,  this integral is bounded above by 
$(1-\cos2\pi\alpha)^{-1/\alpha} \Gamma(1/\alpha)r^{-1/\alpha} $

\end{proof}
{\sl Proof  of Theorem \ref{main}. }
 
Let $G_{\alpha} (-r^+)= R(r)e^{-i\pi\th(r)}$ be the polar decomposition of  $G_{\alpha}(-r^+)$. 
 Since $\Im G_{\al}(-r^+)$ is negative, we can choose $\th(r)\in ]0,1[$ and continuous.  
 Proposition \ref{Galpha} implies that 
\begin{eqnarray*}
\th(r)&=&\alpha +O(r^{1-\alpha})\qquad r\to 0\\
\th(r)&= &1/2 + O(1/r)\qquad r\to+\infty
\end{eqnarray*}
Let 
$$L_\alpha(z)=\exp\int_0^\infty \left[\frac{1}{z+t}-\frac{1}{1+t}\right]\theta(t)dt,$$
this function
 is  analytic  on ${\bf C}\setminus ]-\infty,0]$ and  satisfies,  by well known properties of Stieltjes transforms, 
$$\frac{L_\alpha(-r^+)}{L_\alpha(-r^-)}=e^{-2i\pi\theta(r)}\qquad r>0,$$
furthermore,  since $\theta(t) = 1/2+ O(1/t)$, the integral $\int _0^{\infty} {1\over t+1} (\th(t)- 1/2) dt$  is finite and 
$z^{1/2} L_\alpha(z)=\exp\int_0^{\infty}\left[{1\over z+t} - {1\over 1+t}\right] (\th(t)-1/2) dt$ therefore
$$z^{1/2} L_\alpha(z)\to_{z\to\infty} \exp \int_0^{\infty}{1\over 1+t}  (\th(t)-1/2)dt=C>0$$
and
$L_\alpha(z) \sim  Cz^{-1/2},\quad z\to \infty.$ 
A similar argument, using the fact that $\th(t)=\alpha +O(t^{1-\alpha})\quad t\to 0$, gives 
 $$z^{\alpha}L_\alpha(z)\to _{z\to 0} \exp \int_0^{\infty}  {1\over t(1+t)} (\th(t)-\alpha) dt= D >0$$
On the other hand,
$$\frac{G_\alpha(-r^+)}{G_\alpha(-r^-)}=e^{-2i\pi\theta(r)}$$
therefore the function
$E_\alpha(z)={e^{\delta z}G_\alpha(z)\over L_\alpha(z)}$
is analytic on ${\bf C}\setminus ]-\infty,0]$, and has a  continuous extension to ${\bf C}\setminus\{0\}$.  It is also continuous at $0$, because both $L_{\alpha}(z)$ and $G_{\alpha}(z)$ are equivalent to $z^{-\alpha}$ up to a  multiplicative constant, for $z\to 0$. By Morera's theorem, the function $E_{\alpha}$   can be extended to  an entire function. Moreover, since the two functions $e^{\delta z}G_\alpha(z)$ and $L_\alpha(z)$ are equivalent to $z^{-1/2}$  at infinity up to a multiplicative constant, 
 $E_{\alpha}(z)$ is bounded on {\bf C}. Finally, by Liouville theorem, $E_{\alpha}$ is constant and this constant, equal to $e^{\delta}G_{\alpha} (1)$ is positive. \qed
 
 {\bf Remark} : If  $H$ is an HCM function then ${\log H(x)\over x}$ is bounded, thus 
 $G_{\alpha} (x^h)$ is not HCM for any $h>0$. Consequently,  if $g_{\alpha}(x^{\gamma})$ is HCM  then $\gamma< \beta^{-1}$.

\section{HCM,  non HCM, anti HCM property of $G_{\alpha}$}

\begin{theorem} \label{Thcm}

\begin{enumerate}\label{GHCM}
\item For $\alpha \in ]1/2,1]$,   the function $\th$ is decreasing and $G_{\alpha}$ is anti-HCM, 

\item  For $\alpha \in [1/3,1/2]$, the function $\th$ is increasing and  $G_{\alpha}$ is HCM,   

\item For $\alpha \in ]0,1/3[$, the function $\th$ is not monotonous and $G_{\alpha}$ is neither HCM, neither anti-HCM.
\end{enumerate}
\end{theorem}
For the proof we need some preliminary results.
\begin{lemma}

\begin{enumerate}
\item 
 $r^{\alpha} \Im( G_{\alpha} (-r^+))$ is negative and decreasing 
\item

For $\alpha\in ]1/3,1]$, $sign(1/2-\alpha)r^{\alpha}\Re( G_{\alpha} (-r^+))$ is positive and increasing 
\item
For $\alpha\leq 1/3$, $r^{\alpha^{-1}}\Re(G_{\alpha} (-r^+)$ is not of constant sign. 
\end{enumerate}
\end{lemma}

\begin{proof}
From (\ref{Galt+}) we get: 
$$\Im G_{\alpha}(-r^+))=-(2\pi\beta)^{-1} \int_0^{\infty} e^{-rt} \Re( e^{rt^{\alpha}} - e^{rt^{\alpha} e^{-2i\pi \alpha} })dt $$

$$\Re( G_{\alpha}(-r^+))=(2\pi\beta)^{-1} \int_0^{\infty} e^{-rt} \Im ( e^{rt^{\alpha}} - e^{rt^{\alpha} e^{-2i\pi \alpha} })dt $$

The change of variables $t\to t/r$ in the first identity gives

$$\Im(G_{\alpha}(-r^+))=(2\pi \beta r)^{-1} \int_0^{\infty} \Re(e^{t^{\alpha}r^{1-\alpha}} -   e^{t^{\alpha}r^{1-\alpha}e^{-2i\pi \alpha} })e^{-t }dt$$

while $t\to t/r^{\alpha} $ in the second gives

$$\Re( G_{\alpha}(-r^+))=(2\pi\beta r^{\alpha})^{-1} \int_0^{\infty}  \Im(e^{-r^{1-\alpha\over \alpha} t -e^{i\pi (1-2\alpha)} t^{\alpha} }) dt$$

The function 
$$r^{-1+\alpha}\Re(e^{t^{\alpha}r^{1-\alpha} } -   e^{t^{\alpha}r^{1-\alpha}e^{-2i\pi \alpha}})=\sum_n
\frac{t^{n\alpha}r^{(1-\alpha)(n-1)}}{n!}(1-\cos 2n\pi \alpha)$$ is increasing in $r$ for all $t>0$. It follows that   $r^{\alpha} \Im( G_{\alpha}(-r^+))$ is  increasing.

The second identity and (\ref{Zolotarev}) give

$$r^{\alpha}\Re( G_{\alpha}(-r^+))=\beta^{-1} x^{\alpha+1}g_{\alpha,1-2\alpha} (x)$$ 
( for $x=r^{-{1\over \beta}}$ )

The end of the lemma   follows from  \ref {Monot}.

\end{proof}

{\sl Proof of  theorem \ref{Thcm}} : According to section \ref{HCM} it is enough to consider  monotonicity properties of $\th$. Recall that, for all $\alpha$,  $\th(0)=\alpha$ and $\th(+\infty)=1/2$, moreover
 $\Im (G_{\alpha}(-r^+))= R(r)\sin\pi\th(r)$ is negative and decreasing and, for $\alpha \geq 1/3$, $\Re(G_{\alpha}(-r^+))= R(r)\sin\pi\th(r)$ has constant sign and is monotonous. It follows that  
$G_{\alpha}(-r^+)$ takes all  its values in a  quarter  plane and $\th(r)$ has  a constant sign and its absolute value is increasing, thus $\th$ is monotonous, decreasing  for $\alpha>1/2$ and increasing for $\alpha\in [1/3,1/2]$. 

Finally for $\alpha\in ]0,1/3[$, we obtain that $\Re(G_{\alpha} (-r^+))= R(r)\cos\pi\th(r)$ can take negative values, thus $\th(r)$ does not take all its value inside the interval  $[\alpha, 1/2]$, thus it is not monotonous.

In the case $\alpha \in [1/3,1[$ we obtain a better estimate for $G_{\alpha}(x)$  than in corollary \ref{sim}.

\begin{corollary} Let $$f_1(x)=x^{-\alpha}e^{-\delta x}1_{[0,1]}(x), \quad f_2(x)=x^{-1/2}e^{-\delta x} 1_{]1,+\infty[} $$

If  $\alpha\in [1/3,1/2]$, then

$$ G_{\alpha}(1)f_1(x) +  c_{\infty}f_2(x)
\leq G_{\alpha} (x)\leq    c_0f_1(x)  + G_{\alpha}(1)f_2(x)$$

If $\alpha\in ]1/2,1]$, then 

$$ G_{\alpha}(1)f_1(x) +  c_{\infty}f_2(x)
\geq G_{\alpha} (x)\geq    c_0f_1(x)  + G_{\alpha}(1)f_2(x)$$
\end{corollary}

Using the proposition \ref{HCMGGC}, we also obtain new GGC densities related to $\alpha$- stable densities.

 [$\varepsilon( \alpha-1/2)$ denotes the sign of $\alpha-1/2$].

\begin{corollary} If  $\alpha\in [1/3, 1]$ then

1) The function $[c_0^{-1}.x^{-\alpha}G_{\alpha} (x)]^{-\varepsilon( \alpha-1/2)}$ is the Laplace transform of a random variable  of the form $Y-\varepsilon(\alpha-1/2) \delta$ where $Y$ is GGC.

2) $[c_{\infty}^{-1}x^{-1/2} e^{\delta /x}G_{\alpha} (1/x)]^{\varepsilon( \alpha-1/2)}$ is the Laplace transform of a GGC random variable.

\end{corollary}

Finally we obtain another consequence for the $\alpha$-densities.

\begin{corollary}

If $\alpha>1/2$ then $S_{\alpha}$ is GGC and $S_{\alpha} ^{-\beta}$   is not GGC.

\end{corollary}

{\sl proof } The GGC property of $S_{\alpha}$ is known and has already been already been mentioned  in paragraph 2.2.
Consider the Laplace transform, for   $\lambda\geq 0$, 

$$Lp(S_{\alpha} ^{-\beta})(\lambda) = \int_0^{\infty} e^{-\lambda x} G_{\alpha} (x) dx =\int_0^{\infty} e^{-t} G_{\alpha} ({t\over \lambda} ) {dt\over \lambda}$$

Since $(z^{1/2} 1_{|z|>1} + z^{\alpha} 1_{|z|\leq 1}) e^{\delta z} G_{\alpha} ( z)$ is a bounded function of $z$,  the   integral can be analytically continued by an analytic to $\C\backslash ]-\infty, 0]\cap \{ |z|>\delta \}$ and this continuation satisfies again  for $r>\delta$,

$$ Lp( S_{\alpha} ^{-\beta})(-r^+)=-\int _0^{\infty} e^{-t} G_{\alpha}(-({t \over r})^- )  {dt \over r}=- \int_0^{\infty} e^{-rt} G_{\alpha}  (-t^-)
dt$$

Since $\Im(G_{\alpha} (-t^-))$ is positive  and increasing while $\Re(\Im(G_{\alpha} (-t^-))$ is negative and increasing, the same is true for $- \Im Lp( S_{\alpha} ^{-\beta})(-r^+)$ and $-\Re(Lp( S_{\alpha} ^{-\beta})(-r^+)$, consequently the opposite of the argument of $Lp( S_{\alpha} ^{-\beta})(-r^+)$ is decreasing again for $r>\delta$. 

Thus $Lp( S_{\alpha} ^{-\beta}) (x)$ cannot be HCM, consequently $S_{\alpha} ^{-\beta}$ is not GGC.

\section {Some further properties of  GGC and HCM functions}

 \begin{theorem}\label {HCMGGCb}

Let  $H$ be an   HCM function and  $g$ be a GGC density,  then the function 
$\int_0^{\infty} H(xy) g(y) dy$ is HCM if it is finite. 

\end{theorem}

For the proof of this result we derive some lemmas. The first one is due to Bondesson \cite{bond2}.

\begin{lemma} \label{bondlemm}

The product and  the ratio of  two  independent GGC random variables is GGC.

\end{lemma}

From this we deduce:
\begin{lemma} \label{mom} Let $g$  be a GGC density and $\beta$ a real number such that  $m_{\beta}=\int_0^{\infty} x^{\beta} g(x)dx$ is finite, then 
$m_{\beta}^{-1} x^{\beta} g(x)$ is a GGC density.

\end{lemma} 
\begin{proof}  

Let  $H$ be  the   Laplace transform of a GGC random variable $Y$, and  $X$ be a GGC random variable with density  $g$ independent of $Y$.

The function $\int_0^{\infty} H(xy) g(y) dy$ is the Laplace transform of $XY$. According to lemma \ref{bondlemm} the independent  product $XY$ is GGC again, thus 
$\int_0^{\infty} H(xy) g(y) dy$ is the Laplace transform of a GGC variable. Thus it is HCM. 
Replacing $H(x)$ by  this  to  $H(x)e^{-x} (\varepsilon ^{-\beta} (\varepsilon+x)^{-\beta})$ which is the Laplace  transform of $Y+ E_{\varepsilon}+1$ where $E_{\varepsilon}$ has $\Gamma(\varepsilon, \beta)$-distribution, 
we obtain that  the integral $$\int_0^{\infty} H(xy) \varepsilon ^{-\beta} (\varepsilon+xy)^{-\beta} g(y)dy$$ is HCM. Multiply this integral by the constant $\varepsilon ^{\beta}$ and let  $\varepsilon\to 0$, 
the  monotone convergence  theorem, and the fact that  HCM property  is stable by multiplication by a  positive constant and by  pointwise limit  implies that $$\int_0^{\infty}  e^{-xy} (xy)^{-\beta} g(y)dy$$  is HCM.
Mutiplying  this integral by $m_{\beta}^{-1}x^{\beta}$ we again get an  HCM function and  the integral obtained is the Laplace transform of the density
$m_{\beta}^{-1}x^{-\beta} g(x)$. Since this Laplace transform is HCM the density is GGC.
\end{proof}
 
 \begin{lemma}\label{approx}

Let   $\th$ be an increasing function and $H$ the associated HCM funtion 
$$H(x) =\exp \int_0^{\infty} ({1\over x+t}- {1\over 1+t}) \th(t) dt$$

Then $H$  is a pointwise limit of HCM functions $H_n$  whose  $\th$-function in the  representation (\ref{HCMb})  is bounded.
Moreover one can chose the $H_n$ such that   for all $\epsilon>0$  there exists $N$ s.t.  if  $n>N$ and  $x\in ]0,+\infty[$ then 

$$( 1-\varepsilon)H(x) \leq H_n(x)\leq H(x)e^{\varepsilon (x+x^{-1})}$$ 

\end{lemma}

 \begin{proof}

Let $n$ be a positive integer and
$$\th(t)=\th_n + (\th(t) -n)1_{\th(t)\geq n} + (\th(t) +n)1_{\th(t)\leq -n}$$
with $$\th_n=\th (t)\vee (-n) \wedge n$$

Moreover let 

$$H_n(x)=\exp \int_0^{\infty} ({1\over x+t}- {1\over 1+t}) \th_n(t) dt$$
$$ E_n(x)=\int _0^{\infty} ( {1\over x+t}- {1\over 1+ t}) (\th(t)- n)1_{\th(t)>n}dt $$
$$\hat E_n(x)=\int_0^{\infty}({1\over x+t}- {1\over 1+ t}) (\th(t)+ n)1_{\th(t)<-n}dt $$

Clearly $$H=\hat E_nH_nE_n$$

Since  ${1\over x+t}- {1\over 1+ t} $ and $x-1$ have the same sign
and 
$$ -1_{x>1}\inf(x,t^{-2})\leq {1\over x+t}- {1\over 1+ t} \leq  1_{x<1} \inf (x^{-1},t^{-2})  $$

we obtain

$$e^{-\varepsilon _n x 1_{x>1} }\leq E_n\leq e^{\varepsilon_n 1_{x<1}}$$

and 

$$ e^{-\hat \varepsilon _n x^{-1} 1_{x<1} }\leq \hat E_n\leq e^{\hat \varepsilon_n 1_{x>1}}$$

with  $$\varepsilon_n=\int_0^{\infty} 1_{\th(t)>n}(\th(t)-n){dt\over t^2} $$

and $$ \hat\varepsilon_n= \int_0^{\infty} 1_{\th(t)<-n}(-n-\th(t) )dt$$
The positive numbers  $\varepsilon _n$ and $\hat\varepsilon_n$ go to zero when $n\to+\infty$ for all $\varepsilon >0$ , let $N$ such that $\varepsilon >\varepsilon _n\vee\hat\varepsilon_n$  then $H_n(x)$ satisfies the required estimate of the lemma.
\end{proof}

{\it Proof of Theorem \ref{HCMGGCb}}

Let $X$ and $Y$ be GGC random variables, let $H$ be the 
 Laplace transform of  $Y$ and $g$ be the density of $X$. 
 On the other hand, 
the sequence of functions $(1+{b\over x n})^{-n}= x^n {b\over n}^n ( 1+ {nx\over b})^{-n}$  have limit $e^{-b\over x}$ and they are bounded by 1.
Moreover the function $( 1+ {xn\over b})^{-n}$ is the Laplace transform of a  $\Gamma (n,  {n\over b})$- random variable $E_n$,
The product $H(x)( 1+ {nx\over b})^{-n}$ is the Laplace transform of the GGC variable $Y+E_n$,  since the independant product $X(Y+E_n)$ is again 
GGC and its Laplace transform is 
$$\int_0^{\infty} H(xy) ( 1+ {nxy\over b})^{-n}g(y)dy$$

Thus, this function is HCM. 

Suppose  that the random $X$ has moments of all order, according to lemma \ref{mom}  the function $g(y)$ can be replaced by $g(y)y^{\beta+n}$ for any $n$ and $\beta$, and  again the integral
$$\int_0^{\infty} H(xy) ( 1+ {nxy\over b})^{-n}y^{\beta +n}g(y)dy$$

defines an HCM function. 

Multiply by $x^{n}({b\over n})^n$ and let $n$ goes to infinity and use  the dominated convergence theorem 
in order to obtain that for all real $\beta$ and $b>0$.

$$\int_0^{\infty} H(xy) e^{-b \over xy} y^{\beta} g(y)dy$$ is HCM.

The function $H(x)e^{- ax}$ can replace $H(x)$  in this formula, since it is also a   Laplace transform of  $Y+a$ which is  a GGC variable again and 
$$\int_0^{\infty}H(xy)e^{- a xy}.  e^{-b \over xy} y^{\beta} g(y)dy$$

is also HCM. 

Take  $a>0$ and $b>0$,  the hypothesis that $X$ (with density $g$) has moments of all orders can be removed because the GGC densities with finite moments are dense in the family of GGC densities  for  the weak topology. Finally, the function 

$$\int_0^{\infty}H(xy)e^{- a xy}.  e^{-b \over xy} y^{\beta} g(y)dy$$

is HCM for any $H$ which is the  Laplace transform of a GGC density ,  any real $\beta$ and any GGC-density $g$. 

 Let $H$ be any HCM function of the form 

$$H(x)=\exp \int_0^{\infty} ({1\over x+t}- {1\over 1+t}) \th(t) dt$$

for an increasing function $\th$,  and $(H_n)$ be a sequence of HCM functions approaching $H$ as it is discribed in  lemma \ref{approx} 
 The $\th$ functions of $H_n$  are bounded bellow (say by $-n$), then $H_n$ are of the form $x^n {\tilde H}_n(x)$ where ${\tilde H}_n$ are Laplace  transform of GGC-variables.( see proposition \ref{HCMGGC}).
Thus , the functions 

$$\int_0^{\infty}H_n(xy)e^{- a xy}  e^{-b \over xy} (xy)^n.y^{\beta} g(y)dy$$

are HCM . 

 Divide by $x^{-n}$ and let $\beta=-n$ and  obtain  that the  functions

$$\int_0^{\infty}H_n(xy)e^{- a xy}.  e^{-b \over xy} g(y)dy$$
are HCM. 

Finally, by Lebesgue dominated convergence in $L^1(R^+, e^{- \varepsilon (x+ x^{-1}}))$ with $0<\varepsilon< \inf(a,b)$, the integral $\int_0^{\infty} H_n (xy) e^{- a xy}.  e^{-b \over xy} g(y)dy$ converges for all $x>0$  to $\int_0^{\infty} H(xy)e^{- a xy}.  e^{-b \over xy} g(y)dy$ and this  function is  again HCM.

The monotone convergence theorem enable to extend  the property  in case   $a$ or $b$ are zero and the proof is finished.

\qed

  \begin{corollary} [T-Simon , P.Bosch] 
  
  The $\alpha$ stable positive density, $g_{\alpha}$, is HCM  if and only if $\alpha\in ]0,1/2]$. 
  
  \end{corollary}
  
  \begin{proof} The scaling property gives the identity for $\alpha\leq 1/2$
  
  $S_{\alpha}= S_{1/2}^{1\over 2\alpha}. S_{2\alpha}$. 
  
The density of the standard $1/2$-stable r.v., $S_{1/2}^{1\over 2\alpha}$ is    
  $H(x)=(2\pi)^{-1/2}.\alpha x^{-3\alpha-1} e^{-x^{-2\alpha}}$ which is clearly  HCM .
  
  Applying theorem \ref{HCMGGCb}  to  the HCM function $x^{-1}H(x^{-1})$ and to  the density of  the GGC variable $S_{2\alpha}$   gives the required property.

  \end{proof}
  
 Finally, let  $\gamma$ be the bigger power such that $g_{\alpha}(x^{\gamma})$ is HCM, we have obtained that    $\gamma=\beta^{-1}$ for $\alpha\in [1/3,1/2]$ and $\gamma \in[1,\beta^{-1}[$ for $\alpha<1/3$.


\begin{thebibliography}{10}

\bibitem{bond}
Bondesson, L.{\sl
Generalized gamma convolutions and related classes of distributions and densities.}
Lecture Notes in Statistics, 76. Springer-Verlag, New York, 1992.

\bibitem{bond2} Bondesson,L. {\sl A class of probability distributions that is closed with respect to addition as well as multiplication
of independent random variables}. To appear in Journal of Theoretical Probability, 2015.

\bibitem{Bosch} Bosch, P. Simon, T. {\sl A proof of Bondesson conjecture on stable densities}
Arkiv f\"or Matematik. April 2016, Volume 54, Issue 1, pp 31-38.

\bibitem{Fou} S.Fourati {\sl $\alpha$-stable densities are hyperbolically completely monotone for $\alpha \in ]0,1/4]\cup [1/3,1/2]$}  arXiv:1309.1045 [math.PR].Sept. 2013.

\bibitem{Yor-Roy-James} James, L.F, Roynette B., Yor, M.  {\it Generalized Gamma Convolutions, Dirichlet means, Thorin measures, with explicit examples}.
Probability Surveys
Vol. 5 (2008) 346–415

\bibitem{WissSim}
Jedidi, W. and  Simon T. {\sl Further examples of GGC and HCM functions}. Bernoulli 19 (5A), 1818-1838, 2013.

\bibitem{thor}
Thorin, Olof {\sl On the infinite divisibility of the Pareto distribution.}
Scand. Actuar. J. 1977, no. 1, 31'A ì40.

\bibitem{zol}
Zolotarev, V. M. {\sl
One-dimensional stable distributions.}
Translated from the Russian by H. H. McFaden. Translation edited by Ben Silver. Translations of Mathematical Monographs, 65. American Mathematical Society, Providence, RI, 1986.
  
\end{thebibliography}
\end{document}